


\magnification=1200 

\vsize=18cm \voffset=1cm \hoffset=1cm
\hsize=11.25cm
\parskip 0pt
\parindent=12pt

\catcode'32=9

\font\tenpc=cmcsc10
\font\eightpc=cmcsc8
\font\eightrm=cmr8
\font\eighti=cmmi8
\font\eightsy=cmsy8
\font\eightbf=cmbx8
\font\eighttt=cmtt8
\font\eightit=cmti8
\font\eightsl=cmsl8
\font\sixrm=cmr6
\font\sixi=cmmi6
\font\sixsy=cmsy6
\font\sixbf=cmbx6

\skewchar\eighti='177 \skewchar\sixi='177
\skewchar\eightsy='60 \skewchar\sixsy='60

\font\tengoth=eufm10
\font\tenbboard=msbm10
\font\eightgoth=eufm7 at 8pt
\font\eightbboard=msbm8
\font\sevengoth=eufm7
\font\sevenbboard=msbm7
\font\sixgoth=eufm6
\font\fivegoth=eufm5

\font\tenmsa=msam10
\font\eightmsa=msam8
\font\sevenmsa=msam7

\newfam\gothfam  
\newfam\bboardfam 
\newfam\msafam 

\catcode`@=11

\def\tenpoint{%
  \textfont0=\tenrm \scriptfont0=\sevenrm \scriptscriptfont0=\fiverm
  \def\rm{\fam\z@\tenrm}%
  \textfont1=\teni \scriptfont1=\seveni \scriptscriptfont1=\fivei
  \def\oldstyle{\fam\@ne\teni}%
  \textfont2=\tensy \scriptfont2=\sevensy \scriptscriptfont2=\fivesy
  \textfont\gothfam=\tengoth \scriptfont\gothfam=\sevengoth
  \scriptscriptfont\gothfam=\fivegoth
  \def\goth{\fam\gothfam\tengoth}%
  \textfont\bboardfam=\tenbboard \scriptfont\bboardfam=\sevenbboard
  \scriptscriptfont\bboardfam=\sevenbboard
  \def\bboard{\fam\bboardfam}%
  \textfont\msafam=\tenmsa \scriptfont\msafam=\sevenmsa
  \scriptscriptfont\msafam=\sevenmsa
	\textfont\itfam=\tenit
  \def\it{\fam\itfam\tenit}%
  \textfont\slfam=\tensl
  \def\sl{\fam\slfam\tensl}%
  \textfont\bffam=\tenbf \scriptfont\bffam=\sevenbf
  \scriptscriptfont\bffam=\fivebf
  \def\bf{\fam\bffam\tenbf}%
  \textfont\ttfam=\tentt
  \def\tt{\fam\ttfam\tentt}%
  \abovedisplayskip=12pt plus 3pt minus 9pt
  \abovedisplayshortskip=0pt plus 3pt
  \belowdisplayskip=12pt plus 3pt minus 9pt
  \belowdisplayshortskip=7pt plus 3pt minus 4pt
  \smallskipamount=3pt plus 1pt minus 1pt
  \medskipamount=6pt plus 2pt minus 2pt
  \bigskipamount=12pt plus 4pt minus 4pt
  \normalbaselineskip=12pt
  \setbox\strutbox=\hbox{\vrule height8.5pt depth3.5pt width0pt}%
  \let\bigf@ntpc=\tenrm \let\smallf@ntpc=\sevenrm
  \let\petcap=\tenpc
  \normalbaselines\rm}

\def\eightpoint{%
  \textfont0=\eightrm \scriptfont0=\sixrm \scriptscriptfont0=\fiverm
  \def\rm{\fam\z@\eightrm}%
  \textfont1=\eighti \scriptfont1=\sixi \scriptscriptfont1=\fivei
  \def\oldstyle{\fam\@ne\eighti}%
  \textfont2=\eightsy \scriptfont2=\sixsy \scriptscriptfont2=\fivesy
  \textfont\gothfam=\eightgoth \scriptfont\gothfam=\sixgoth
  \scriptscriptfont\gothfam=\fivegoth
  \def\goth{\fam\gothfam\eightgoth}%
  \textfont\bboardfam=\eightbboard \scriptfont\bboardfam=\sevenbboard
  \scriptscriptfont\bboardfam=\sevenbboard
  \def\bboard{\fam\bboardfam}%
  \textfont\msafam=\eightmsa \scriptfont\msafam=\sevenmsa
  \scriptscriptfont\msafam=\sevenmsa
	\textfont\itfam=\eightit
  \def\it{\fam\itfam\eightit}%
  \textfont\slfam=\eightsl
  \def\sl{\fam\slfam\eightsl}%
  \textfont\bffam=\eightbf \scriptfont\bffam=\sixbf
  \scriptscriptfont\bffam=\fivebf
  \def\bf{\fam\bffam\eightbf}%
  \textfont\ttfam=\eighttt
  \def\tt{\fam\ttfam\eighttt}%
  \abovedisplayskip=9pt plus 2pt minus 6pt
  \abovedisplayshortskip=0pt plus 2pt
  \belowdisplayskip=9pt plus 2pt minus 6pt
  \belowdisplayshortskip=5pt plus 2pt minus 3pt
  \smallskipamount=2pt plus 1pt minus 1pt
  \medskipamount=4pt plus 2pt minus 1pt
  \bigskipamount=9pt plus 3pt minus 3pt
  \normalbaselineskip=9pt
  \setbox\strutbox=\hbox{\vrule height7pt depth2pt width0pt}%
  \let\bigf@ntpc=\eightrm \let\smallf@ntpc=\sixrm
  \let\petcap=\eightpc
  \normalbaselines\rm}

\def\pc#1#2|{{\bigf@ntpc #1\penalty \@MM\hskip\z@skip\smallf@ntpc%
    \uppercase{#2}}}

\catcode`@=12

\tenpoint

\def\pointir{\discretionary{.}{}{.\kern.35em---\kern.7em}\nobreak
   \hskip 0em plus .3em minus .4em }

\def\qed{\quad\raise -2pt\hbox{\vrule\vbox to 10pt{\hrule width 4pt
   \vfill\hrule}\vrule}}

\def\rem#1|{\par\medskip{{\it #1}\pointir}}

\def\abstract#1{\vbox{\eightpoint\narrower\narrower
\pc ABSTRACT|.\quad #1}}

\def\section#1{\goodbreak\par\vskip .3cm\centerline{\bf #1}
   \par\nobreak\vskip 3pt }

\def\article#1|#2|#3|#4|#5|#6|#7|
    {{\leftskip=7mm\noindent
     \hangindent=2mm\hangafter=1
     \llap{{\tt [#1]}\hskip.35em}{\petcap#2}\pointir
     #3, {\sl #4}, {\bf #5} ({\oldstyle #6}),
     pp.\nobreak\ #7.\par}}

\def\livre#1|#2|#3|#4|
    {{\leftskip=7mm\noindent
    \hangindent=2mm\hangafter=1
    \llap{{\tt [#1]}\hskip.35em}{\petcap#2}\pointir
    {\sl #3}, #4.\par}}

\def\divers#1|#2|#3|
    {{\leftskip=7mm\noindent
    \hangindent=2mm\hangafter=1
     \llap{{\tt [#1]}\hskip.35em}{\petcap#2}\pointir
     #3.\par}}


\def\Sym{{\goth S}}

\mathchardef\Cup="2A64  
\def\des{\mathop{\rm des}\nolimits}

\def\des{\mathop{\rm des}\nolimits}
\def\maj{\mathop{\rm maj}\nolimits}
\def\ides{\mathop{\rm ides}\nolimits}
\def\inv{\mathop{\rm inv}\nolimits}
\def\imaj{\mathop{\rm imaj}\nolimits}

\def\Mc{\mathop{\rm Mc}\nolimits}
\def\McInv{\mathop{\rm Mc}^{-1}\nolimits}
\def\IcInv{\mathop{\rm Ic}^{-1}\nolimits}
\def\Sc{\mathop{\rm Sc}\nolimits}
\def\Lc{\mathop{\rm Lc}\nolimits}
\def\Ic{\mathop{\rm Ic}\nolimits}
\def\Ligne{\mathop{\rm Ligne}\nolimits}
\def\Iligne{\mathop{\rm Iligne}\nolimits}
\def\Majcode{\mathop{\rm Majcode}\nolimits}
\def\Invcode{\mathop{\rm Invcode}\nolimits}
\def\Eul{\mathop{\rm Eul}\nolimits}
\def\El{\mathop{\rm El}\nolimits}
\def\eul{\mathop{\rm eul}\nolimits}
\def\sort{\mathop{\rm sort}\nolimits}

\def\SE{\mathop{{\hbox{SE}}}\nolimits}


\def\last{\mathop{\hbox{\rm last}}}
\def\first{\mathop{\hbox{\rm first}}}
\def\id{\mathop{\hbox{\rm id}}\nolimits}

\def\bfr{\mathop{\bf r}}
\def\bfi{\hbox{\bf i}\thinspace}
\def\bfc{\mathop{\bf c}}
\def\compose{\circ}
\def\Longmapsto#1{\hbox to #1{$\mapstochar$\rightarrowfill}}

%


\def\RefDW{DW93}
\def\RefFo{Fo68}
\def\RefFH{FH04}
\def\RefFS{FS78}
\def\RefGG{GG79}
\def\RefHa{Ha90}
\def\RefHNT{HNT06}
\def\RefLe{Le60}
\def\RefLoI{Lo83}
\def\RefLoII{Lo02}
\def\RefMac{Mac15}
\def\RefRa{Ra81}
\def\RefSk{Sk01}
\def\ThCompl{2}
\def\ThFactIncr{3}
\def\ThCycleMc{4}
\def\ThSortMc{5}
\def\ThElRec{6}
\def\ThShuffle{7}
\def\ThMajcodeMc{8}
\def\ThElEulA{9}
\def\ThElEulB{10}
\def\SecCompl{3}

\rightline{2006/12/22}

\bigskip \bigskip \bigskip
\centerline{\bf Euler-Mahonian triple set-valued statistics on permutations}
\bigskip 
\centerline{Guo-Niu Han}
\medskip
\centerline{Center for Combinatorics, LPMC}
\centerline{Nankai University, Tianjin 300071, P. R. China}
\medskip
\centerline{I.R.M.A. UMR 7501, Universit\'e Louis Pasteur et CNRS}
\centerline{7, rue Ren\'e-Descartes, F-67084 Strasbourg, France}
\smallskip
\centerline{\tt guoniu@math.u-strasbg.fr}

\bigskip \bigskip
\abstract{
The inversion number and the major index are equidistributed on
the symmetric group. This is a classical result, first proved by 
MacMahon [\RefMac],	then by Foata by means of a combinatorial 
bijection [\RefFo]. Ever since many refinements have been derived, 
which consist of adding new statistics, or replacing integral-valued statistics by set-valued ones. 
See the works by Foata-Sch\"utzen\-ber\-ger [\RefFS],  Skandera [\RefSk], 
Foata-Han [\RefFH] and  more recently by Hivert-Novelli-Thibon [\RefHNT].
In the present paper we derive a general equidistribution property on  
Euler-Mahonian set-valued statistics on permutations, which unifies 
the above four refinements.
We also state and prove the so-called ``complement property" of the Majcode.
}
\medskip

\bigskip
\section{1. Introduction} 
Let $w=y_1y_2\cdots y_n$ be a word whose letters $y_1, y_2, \ldots, y_n$
are integers. 
The {\it descent number} ``des",  {\it major index} ``maj" and 
{\it inversion number} ``inv" are defined by 
(see, for example, [Lo83, \S 10.6]):
$$
\eqalign{
\des w&=\#\{i \mid 1\leq i\leq n-1, y_i>y_{i+1} \},\cr
\maj w&=\sum\{i \mid 1\leq i\leq n-1, y_i>y_{i+1} \},\cr
\inv w&=\#\{(i,j) \mid 1\leq i<j\leq n, y_i>y_j \}.\cr
}
$$
In this paper we only deal with permutations $\sigma=x_1x_2\cdots x_n$
of $12\cdots n$ $(n\geq 1)$. 
A statistic is said to be {\it Mahonian}, if it has the same distribution as 
``maj" on the symmetric group $\Sym_n$,
and 
a bi-statistic is said to be {\it Euler-Mahonian} if it has the same distribution as 
(des, maj). 
MacMahon's 
fundamental result says that ``inv" is Mahonian [\RefMac], i.e.,
``maj" and ``inv" have the same distribution on  $\Sym_n$.
This equidistribution property will be written
$$
\maj\simeq\inv,\leqno{(M1)}
$$
which also  means that we have:
$$
\sum_{\sigma\in\Sym_n}q^{\maj\sigma}
=
\sum_{\sigma\in\Sym_n}q^{\inv\sigma}.
$$
Foata [\RefFo] obtained a combinatorial proof of MacMahon's result by
constructing an explicit transformation $\Phi$
such that $\maj\sigma=\inv\Phi(\sigma)$.
Let the {\it ligne of route}
of a permutation $\sigma=x_1x_2\cdots x_n$ be the set of all descent 
places:  
$$\Ligne\sigma=\{i\mid 1\leq i\leq n-1, x_i>x_{i+1}\}.$$
The {\it inverse ligne of route} of $\sigma$ is defined by  
$\Iligne\sigma=\Ligne\sigma^{-1}$.
Foata and Sch\"utzenberger [\RefFS] showed that the  
transformation $\Phi$ preserved the inverse ligne of route and then derived
the first refinement of MacMahon's result: 
$$
(\Iligne, \maj)\simeq(\Iligne, \inv).\leqno{(M2)}
$$
A word $w=d_1d_2\cdots d_n$ is said to be {\it subexcedent} if 
$0\leq d_i\leq i-1$ for all $i=1,2,\ldots, n$. The set of all
subexcedent words of length $n$ is denoted by~$\SE_n$.
The Lehmer code [\RefLe] is a bijection 
$\Invcode:\Sym_n\rightarrow\SE_n$  which
maps each permutation $\sigma=x_1x_2\cdots x_n$ onto a
subexcedent word $\Invcode\sigma=d_1d_2\cdots d_n$, where $d_i$
is given by
$$
d_i=\#\{ j\mid 1\leq j\leq i-1, x_j>x_i\}.
$$
The {\it major index code}, denoted by
$\Majcode$, is a bijection of $\Sym_n$ onto $\SE_n$,
which maps each permutation $\sigma=x_1x_2\cdots x_n$ onto a
subexcedent word $\Majcode\sigma=d_1d_2\cdots d_n$, where $d_i$
is given by
$$
d_i=\maj (\sigma|_{i}) -\maj (\sigma|_{i-1}).
$$
In the above expression $\sigma|_{i}\in\Sym_i$ is the permutation
derived from $\sigma$ by erasing the letters  $i+1, i+2, \ldots, n$.
For example, $\Majcode(175389642)=002135573$ 
and $\Invcode(784269135)=002320654$.

Furthermore, ``eul" is an integral-valued statistic (see [\RefHa, \RefFH])
defined on $\SE_n$ as follows. Let
$w=d_1d_2\cdots d_n$ be a subexcedent word. If
$n=1$, then $\eul w=0$; if $n\geq 2$ let 
$w'=d_1d_2\cdots d_{n-1}$ so that $w=w'd_n$, then define
$$
\eul(w)=\cases{
\eul w', &if $d_n\leq \eul w'$, \cr
1+\eul w', &if $d_n\geq 1+ \eul w'$. \cr
}
$$
Skandera [\RefSk] proved the following refinement
$$
(\des,\maj) \simeq(\eul\compose \Invcode, \inv).\leqno{(M3)}
$$
He also conjectured the following multi-variable equidistribution: 
$$
(\des,\maj, \ides, \imaj) \simeq(\des, \maj, \eul\compose \Invcode, \inv),
\leqno{(M4)}
$$
where $\ides\sigma=\des\sigma^{-1}=\#\Iligne\sigma$ 
and $\imaj\sigma=\maj\sigma^{-1}=\sum\Iligne\sigma$.
This conjecture was proved by Foata and Han [\RefFH]. In fact,
we have obtained the following stronger refinement
$$
(\Iligne, \Eul\compose\Majcode) \simeq(\Ligne, \Eul\compose\Invcode),
\leqno{(M5)}
$$
where ``Eul" is a set-valued statistic defined for each subexcedent word,
having the property: 
$\#\Eul=\eul$.  The explicit definition of ``Eul" can be found in [\RefFH].
We also have the alternate definition:
$$
\Ligne\sigma=\Eul\compose\Majcode\sigma.
$$     
Note that there is no ``perfect" {\it vector-based} refinement of MacMahon's 
result because
$$
(\Iligne, \Majcode) \not\simeq(\Ligne, \Invcode).
$$
We only have the {\it set-based} equidistribution displayed in $(M5)$.
\medskip

Recently, another set-based refinement of MacMahon's result 
was  discovered  by Hivert, Novelli and Thibon [\RefHNT].
Their notations are slightly different: they use {\it subdiagonal}
instead of {\it subexcedent} words.
A word $w=d_1d_2\cdots d_n$ is said to be {\it subdiagonal}, if 
$0\leq d_i\leq n-i$ for all $i=1,2,\ldots, n$. 
Instead of ``Invcode" they introduce the ``Lc-code", 
denoted by ``Lc", which
is a bijection that 
maps each permutation $\sigma=x_1x_2\cdots x_n$ onto a
subdiagonal word $\Lc\sigma=d_1d_2\cdots d_n$, where $d_i$
is given by
$$
d_i=\#\{ j\mid i+1\leq j\leq n, x_i>x_j\}.
$$
Let $\Ic\sigma=\Lc(\sigma^{-1})$.
Their
variation of ``Majcode", called ``Mc-code",
denoted by
``Mc", is a bijection that 
maps each permutation $\sigma=x_1x_2\cdots x_n$ onto a
subdiagonal word $\Mc\sigma=d_1d_2\cdots d_n$, where $d_i$
is given by
$$
d_i=\maj (\sigma^{(i)}) -\maj (\sigma^{(i+1)}).
$$
In the above expression $\sigma|^{(i)}$ is the subword of 
$\sigma$ obtained by erasing the letters {\it smaller} than $i$. 
The relations between ``Invcode" and ``Ic" (resp. between 
``Majcode" and ``Mc") are given in  Section 3. 

For each word $w$ let ``$\sort w$" be the nondecreasing rearrangement of~$w$.
Then the result obtained by Hivert {\it et al.} [\RefHNT] is a set-based equidistribution property, which can be rephrased as
$$
(\Iligne, \sort\compose\Mc)\simeq(\Iligne, \sort\compose\Ic).
\leqno{(M6)}
$$

The variation of ``Eul" is denoted by ``El". In this paper 
we simply define ``El" by
$$
\Ligne\sigma=\El\compose\Mc\sigma.
$$     
Some relations between the statistics ``El" and  ``Eul" 
are given in Section~3.

The main result of the present paper is
the following set-based 
equidistribution property, which includes all previous equidistribution
properties
$(M1)-(M6)$ as special cases.  
\proclaim Theorem 1.
The following
two triplets of set-valued statistics are equidistributed
on the symmetric group~$\Sym_n$:
$$
(\Iligne,\sort\compose\Mc, \El\compose\Mc)\simeq
(\Iligne,\sort\compose\Ic, \El\compose\Ic).
\leqno{(M7)}
$$ 

{\it Remark}. Theorem 1 is not an automatic consequence of $(M6)$.
For example, as shown in [\RefHNT], there is another statistic called ``Sc",
 which also
satisfies 
$$
(\Iligne, \sort\compose\Mc)\simeq(\Iligne, \sort\compose\Sc),
$$
but
$$(\Iligne,\sort\compose\Mc, \El\compose\Mc)\not\simeq
(\Iligne,\sort\compose\Sc, \El\compose\Sc).$$ 

Theorem  1 is proved in Section 2. 
To illustrate the above equidistributions we have listed the twenty-four 
permutations of order 4 and their corresponding 
statistics in Section~4.

\medskip

We also use the transformations $\bfi$, $\bfc$ and $\bfr$ of the
dihedral group. Recall that $\bfi \sigma = \sigma^{-1}$ is the inverse of
the permutation~$\sigma$; then,
$\bfc$ is the {\it complement to $(n+1)$} and $\bfr$ the
{\it reverse image} which map each permutation
$\sigma$, written as a linear word $\sigma=x_1\ldots
x_n$, onto
$$\leqalignno{
\bfc\,\sigma&:=(n+1-x_1)(n+1-x_2)\ldots
(n+1-x_n),\cr
\bfr\,\sigma&:=x_n\,\ldots\,x_2x_1,\cr}
$$
respectively.
For each
subexcedent word 
$w=d_1d_2\cdots d_n\in\SE_n$ let 
$$
\delta w=(0-d_1)(1-d_2)(2-d_3)\cdots (n-1-d_n).
$$
Clearly $\delta w$ also belongs to $\SE_n$. The map $\delta$
is called the {\it complement} map.

The second result of this paper is the ``complement" property
of the Majcode.  
As well-known,
the generating polynomial for the major index over the symmetric group
$\Sym_n$ is equal to 
$$F(q)=(1+q)(1+q+q^2)\cdots (1+q+q^2+\cdots   + q^{n-1}),$$
a polynomial having
symmetric coefficients: $F(q)=q^{n(n-1)/2}F(1/q)$. 
This fact can be checked by constructing a bijection $\sigma \mapsto \tau$
satisfying the relation 
$$\maj\sigma={n(n-1)\over 2} - \maj\tau.\leqno{(R1)} $$
This value-based relation has {\it two} array-based  refinements. As 
the major index is equal to the sum of all descent positions, the
following relation implies $(R1)$~:
$$
\Ligne\sigma=\{1,2,\ldots, n-1\} \setminus  \Ligne\tau.  \leqno{(R2)}
$$
The major index is also the sum of all elements in the Majcode. Therefore
the following relation also implies $(R1)$~:
$$
\Majcode\sigma=\delta\compose \Majcode\tau. \leqno{(R3)}
$$
Classically, there is a trivial way of defining a bijection 
satisfying $(R2)$:  just take $\tau=\bfc\sigma$.
But this way provides no simple relation between $\Majcode\sigma$ 
and $\Majcode\bfc\sigma$.
The following result shows that the complement of $\Majcode$ $(R3)$ is  stronger
that the complement of Ligne~$(R2)$. 

\proclaim Theorem \ThCompl.
Let $\sigma$ and $\tau$ be two permutations satisfying relation $(R3)$.
Then relation $(R2)$ holds.

In Section 3 we first give 
an example 
serving to illustrate
the complement property of Majcode. 
Then we prove Theorem \ThCompl.
Finally, we show how to derive $(M5)$ from Theorems~1 and 2.

\section{2. Proof of Theorem 1}  

The basic idea of the proof is to use the inclusion-exclusion principle, as in
[\RefHNT] or in [\RefDW]. 
We begin with some technical lemmas about the Mc-code.

\proclaim Lemma \ThFactIncr. 
Let $\sigma=x_1x_2\cdots x_n$ be a permutation and 
let $\Mc(\sigma)=d_1d_2\cdots d_n$ be its Mc-code.
If $d_1d_2\cdots d_k$ is nondecreasing for some integer $k$
such that $1\leq k\leq n$,
then
all factors of $\sigma$,
whose letters are less than or equal to $k$, are increasing.

{\it Example}. Take $\sigma=12\ 5\ 9\ 6\ 13\ 3\ 4\ 8\ 1\ 2\ 7\ 10\ 11$ and
$k=7$,
we have $\Mc(\sigma)=0011344020200$. The word $d_1d_2\cdots d_k=0011344$
is nondecreasing. There are four maximal factors 
whose letters are in $\{1,2,\ldots,7\}$:
``5", ``6", ``3~4" and ``1~2~7". They are all increasing.
\medskip

{\it Proof}.
Call {\it bad pair} of $\sigma$ each pair $(y,z)$ of letters 
such that
$1\leq y<z\leq k$ and such that
$zwy$ is a {\it factor} of $\sigma$ 
having all its letters in  $\{1,2,\ldots,k\}$.
It suffices to prove that there is no bad pair in the permutation~$\sigma$.
If $\sigma$ contains some bad pairs, 
let $(y,z)$ be the {\it maximal} bad pair, which means
that $(y',z')$ is not a bad pair for every $y'>y$ or $y'=y$ and $z'>z$. 
Consider the permutation $\sigma'$ obtained from $\sigma$ by deleting
all letters smaller than $y$. 
Then $\Mc(\sigma')=d_yd_{y+1}\cdots d_z\cdots d_n$.
This means that $(y,z)$ is also a bad pair of $\sigma'$.
In fact, all bad pairs of $\sigma'$ are of form $(y, \cdot)$
because $(y,z)$ is the maximal bad pair of $\sigma$.

Let $\sigma'=azyb$ with $a,b$ two factors of $\sigma'$ (we can check that
$z$ is just on the left of $y$ because $(y,z)$ is the maximal bad pair 
of $\sigma$).
When making 
the computation of $\Mc(\sigma')$, consider the insertion 
of the letter $z$. 
Let $a'$ (resp. $b'$) the subword obtained from $a$ (resp. $b$) by deleting
all letters smaller than $z$. 
Then
$$
\eqalign{
d_z&=\maj(a'zb')-\maj(a'b')\cr
&=\cases{
\des(b') &if $\last(a')>\first(b')$ or $|a'|=0$; \cr
|a'|+\des(b') &if $\last(a')<\first(b')$ or $|b'|=0$. \cr
}\cr
}
$$ 
In the above equation $\first(w)$ (resp. $\last(w)$)
denotes the first (or leftmost) (resp. last (or rightmost)) letter of $w$.
In the same manner, we have
$$
\eqalign{
d_y&=\maj(azyb)-\maj(azb)\cr
&=\cases{
\des(b) &if $z>\first(b)$ ; \cr
|az|+\des(b) &if $z<\first(b)$ or $|b|=0$. \cr
} \cr
}
$$ 
However $z>\first(b)$ is not possible, because $(y,z)$ is the maximal bad pair
of $\sigma$.
Hence,
$d_y=|az|+\des(b)\geq |a'|+1+\des(b')>d_z$,
a contradiction.~\qed

\proclaim Lemma \ThCycleMc. 
Let $\sigma$ be a permutation and 
$\Mc(\sigma)=d_1d_2\cdots d_n$.
Let $k\in\{1,2,\ldots,n\}$ be an integer satisfying the following
conditions: 
\smallskip\noindent
\item{}\hbox{\rm (C1)} $d_1\leq d_2 \leq \cdots \leq d_{k-1}$ ;
\item{}\hbox{\rm (C2)} $d_{k-1} > d_k$ ; 
\item{}\hbox{\rm (C3)} $d_k\leq d_1$ ;  
\item{}\hbox{\rm (C4)} $k$ is on the right of all $i<k$ in $\sigma$.
\smallskip\noindent
Then 
$$\Ligne\sigma=\Ligne\tau,$$
where
$\tau=\McInv(d_kd_1d_2\cdots d_{k-1} d_{k+1}d_{k+2}\cdots d_n)$.

{\it Proof}. In fact, $\tau$ can be constructed by means of 
an explicit algorithm.
First, define $\tau'$ by the following steps:

(T1) $\tau'(i)=\sigma(i)$, if $\sigma(i)\geq k+1$;

(T2) $\tau'(i)=\sigma(i)+1$, if $\sigma(i)\leq k-1$;

(T3) $\tau'(i)=1$, if $\sigma(i)=k$;

Then the permutation $\tau$ is obtained from $\tau'$ by making the 
following modifications:

(T4) rearrange the maximal factor of $\tau'$ containing ``1" and 
having all its letters in $\{1,2,\ldots, k\}$
in increasing order.
\medskip

{\it Example}.
Take $\sigma=5\ 6\ 12\ 4\ 10\ 2\ 3\ 9\ 11\ 1\ 7\ 8$ and $k=7$, we have
$\Mc(\sigma)=011233{\bf 0}{}42010$. The following calculation shows
that
$\tau=6\ 7\ 12\ 5\ 10\ 3\ 4\ 9\ 11\ 1\ 2\ 8$. We have
$\Mc(\tau)={\bf 0}011233{}42010$.

$$
\matrix{
\sigma=       & 5& 6& 12& 4& 10& 2& 3& 9& 11& 1& 7& 8& \cr
\noalign{\smallskip}
\noalign{\hrule}
\noalign{\smallskip}
\hbox{(T1)}   &  &  & 12&  & 10&  &  & 9& 11&  &  & 8& \cr
\hbox{(T2)}   & 6& 7& 12& 5& 10& 3& 4& 9& 11& 2&  & 8& \cr
\hbox{(T3)}   & 6& 7& 12& 5& 10& 3& 4& 9& 11& 2& 1& 8& \cr
\hbox{(T4)}   & 6& 7& 12& 5& 10& 3& 4& 9& 11& 1& 2& 8&=\tau\cr
}
$$
\medskip

All factors of $\sigma$ and $\tau$ having their 
letters in $\{1,2,\ldots,k\}$ are increasing, thanks to Lemma 
\ThFactIncr\ and condition (C4).
Therefore, $\Ligne\sigma=\Ligne\tau$ by~(T1).

Let
$\Mc(\tau)= f_1f_2\cdots f_k d_{k+1}d_{k+2}\cdots d_n$.
We need prove
(N1) $f_1=d_k$ and
(N2) $f_{i+1}=d_i$ for $i=1,2,\ldots, k-1$.
We first prove the following property related to the insertion 
of $k$ in $\sigma$.

\smallskip
{(P1). \it Let $\sigma'$ be the word obtained from $\sigma$ by deleting 
all letters
smaller than $k$. Then $\sigma'=\cdots xky \cdots$ or $\sigma'=ky\cdots $
with $x>y>k$.}
\smallskip

{\it Proof of \hbox{\rm (P1)}}.
If $k$ is not the last letter of $\sigma'$, i.e., $\sigma'=\cdots ky\cdots$,
then $y>k$ by (C4). We need prove that 
$\sigma'=\cdots xky\cdots$ (with $x<y$) and $\sigma'=\cdots xk$ are not 
possible. 
If those cases occur,
consider the insertion of $k-1$ into $\sigma'$:
$(k-1)$ is on the left of $k$ by (C4), 
so that $d_{k-1}\leq d_k$; a contradiction with (C2).\qed
\smallskip

{\it Proof of \hbox{\rm (N1)}}.
By Property (P1) and Lemma \ThFactIncr\ 
the permutation $\sigma$ must have the form
$\sigma=azx_1x_2\cdots x_rkyb$ with
$z>y>k>x_r>\cdots x_2>x_1$ ($az$ and $b$ being possibly empty)
and the factor $b$ has all its letters greater than~$k$
because of condition (C4).
Then
$\tau=uz1(x_1+1)(x_2+1)\cdots (x_r+1)yb$ 
by definition of~$\tau$. 
Let $a'$ be the word obtained from $a$ by deleting all letters
smaller than $k$. Then 
$$
d_k=\maj(a'zkyb)-\maj(a'zyb)=\des(yb).
$$
On the other hand,
$$
\eqalign{
f_1&=\maj(\tau)
-\maj(uz(x_1+1)(x_2+1)\cdots (x_r+1)yb)\cr
&=\des((x_1+1)(x_2+1)\cdots (x_r+1)yb)\cr
&=\des(yb)=d_k.\qed\cr
}
$$

{\it Proof of \hbox{\rm (N2)}}.
For $i=1,2,\ldots,k-1$ let $\tau=a(i+1)b$.
Let $\overline a$ (resp. $\overline b$) be the word obtained from 
$a$ (resp. $b$) by deleting all letters smaller than $i+1$.
Let $\hat a$ (resp. $\hat b$) be the word obtained from 
$\overline a$ (resp. $\overline b$) by replacing~$j$ 
by $j-1$ for $j\leq k$. Note that $1\not\in\overline a(i+1)\overline b$
and $k\not\in\hat a i \hat b$. Then 
$$
\eqalign{
f_{i+1}&=\maj(\overline a(i+1)\overline b)-\maj(\overline a\overline b)\cr
&=\maj(\hat a i\hat b)-\maj(\hat a \hat b)\cr
&=\cases{
\des(\hat b) &if $\last(\hat a)>\first(\hat b)$ or $|\hat a|=0$; \cr
|\hat a|+\des(\hat b) &if $\last(\hat a)<\first(\hat b)$ or $|\hat b|=0$. \cr
}\cr
}
$$
Let $\sigma=uiv$. Let $u'$ (resp. $v'$) be the word
obtained from $u$ (resp. $v$) by deleting all letters smaller than $i$.
Note that $k\in u'v'$. Then 
$$
\eqalign{
d_i&=\maj(u'iv')-\maj(u'v')\cr
&= \cases{
\des(v') &if $\last(u')>\first(v')$ or $|u'|=0$; \cr
|u'|+\des(v') &if $\last(u')<\first(v')$ or $|v'|=0$. \cr
}\cr
}
$$
In fact, by definition of $\tau$, we have $\hat a=u'$.
We verify that
$\hat b$ is the word obtained from~$v'$ by removing the letter $k$. 
By Property (P1) and condition~(C4), we have only the following cases:
$v'=\cdots xky\cdots$ (with $x>y>k$),  $v'=\cdots xky\cdots$ (with $x<k<y$),
$v'=ky\cdots$ (with $k<y$) and $v'=\cdots xk$ (with $x<k$). In all those
cases,
$$
\cases{
|\hat a|=|u'|; \cr
\des(\hat b)=\des(v').\cr
}
$$
If $|\hat a|=|u'|=0$, then $f_{i+1}=d_i$. 
If $\first(v')\not=k$, then $\first(v')=\first(\hat b)$ and $f_{i+1}=d_i$.
If $\hat a=u'=\cdots x$ and $v'=ky\cdots$,
then
$(x>k) \Leftrightarrow (x>y)$
by Property (P1). 
Hence, $f_{i+1}=d_i$.
\qed
\smallskip
This ends the proof of Lemma \ThCycleMc.\qed

\proclaim Lemma \ThSortMc. 
Let $\beta$ be a permutation of $\{k+1, k+2, \ldots, n\}$ and let $\sigma$
be a shuffle of $12\cdots k$ and $\beta$ whose Mc-code reads
$$\Mc(\sigma)=d_1d_2\cdots d_kd_{k+1}d_{k+2}\cdots d_n.$$
Then 
$$\Ligne\sigma=\Ligne\tau,$$
where
$\tau=\McInv(\sort(d_1d_2\cdots d_k) d_{k+1}d_{k+2}\cdots d_n)$.

{\it Proof}. By induction.
Define
$$\tau_i=\McInv(\sort(d_1d_2\cdots d_i) d_{i+1}\cdots d_kd_{k+1}d_{k+2}\cdots d_n),$$
so that $\tau_1=\sigma$ and $\tau_k=\tau$.
By definition of $\Mc$, 
we have 
$$d_i \geq \max\{d_1, d_2, \ldots, d_{i-1}\} \hbox{\rm \ or\ }
d_i \leq \min\{d_1, d_2, \ldots, d_{i-1}\}$$
for every $i\leq k$,
because $k$ is on the right
of all letters smaller than $k$
(see also the proof of 
Lemma 6.5 in [\RefHNT]). In both cases  $\Ligne\tau_{i-1}=\Ligne\tau_i$ 
for $2\leq i\leq k$ by Lemma \ThCycleMc. \qed
\smallskip

{\it Example}.
Take $n=12, k=7, \beta=12\ 10\ 9\ 11\ 8$. Let $\sigma$ be the
following shuffle of $1234567$ and $\beta$:
$$\sigma=12\ 1\ 2\ 3\ 10\ 4\ 9\ 5\ 6\ 11\ 7\ 8.$$
Then $\Mc(\sigma)=333214042010$.
The following calculation shows that
$\tau=\tau_7=12\ 4\ 5\ 6\ 10\ 3\ 9\ 2\ 7\ 11\ 1\ 8$. 
$$
\matrix{
\tau_1 =\tau_2=\tau_3
&= \McInv(333214042010) &=&
12\ 1\ 2\ 3\ 10\ 4\ 9\ 5\ 6\ 11\ 7\ 8 \cr
\tau_4 
&= \McInv(233314042010) &=&
12\ 2\ 3\ 4\ 10\ 1\ 9\ 5\ 6\ 11\ 7\ 8 \cr
\tau_5 =\tau_6 
&= \McInv(123334042010) &=&
12\ 3\ 4\ 5\ 10\ 2\ 9\ 1\ 6\ 11\ 7\ 8 \cr
\tau_7
&= \McInv(012333442010) &=&
12\ 4\ 5\ 6\ 10\ 3\ 9\ 2\ 7\ 11\ 1\ 8 \cr
}
$$
We check that
$\Ligne(\sigma)=\Ligne(\tau)=\{1,5,7,10\}$.
\medskip


\proclaim Lemma \ThElRec. 
Let $a,b, c$ be words such that $a,b,ca,cb$ are subdiagonal.
If $\El(a)=\El(b)$, then
$$
\El(ca)=\El(cb).
$$

{\it Proof}. 
By induction we need only prove the lemma when  $c=x$ is a one-letter 
word.
Let $\sigma$ (resp. $\tau$) be the permutation such that $\Mc(\sigma)=xa$
(resp. $\Mc(\tau)=xb$).
Also let $\sigma^{(2)}$ (resp. $\tau^{(2)}$) be the subword obtained from  
$\sigma$ (resp. from $\tau$) by erasing the letter $1$. Then 
$\El(a)=\El(b)$
implies $\Ligne(\sigma^{(2)})=\Ligne(\tau^{(2)})$. 
Now $x=\maj(\sigma)-\maj(\sigma^{(2)}) = \maj(\tau) -\maj(\tau^{(2)})$,
so that the letter $1$ is inserted into $\sigma^{(2)}$ and 
$\tau^{(2)}$ at the same position. Hence, $\Ligne\sigma=\Ligne\tau$.
\qed
\medskip

Let $\alpha\in\Sym_k$ and $\beta=y_1y_2\cdots y_\ell\in\Sym_\ell$ be two permutations. 
A permutation $\sigma\in\Sym_{k+\ell}$ is said to be a {\it shifted shuffle} 
of $\alpha$ and $\beta$, if 
the subword of~$\sigma$, whose letters are $1,2,\ldots,k$ (resp.
$k+1, k+2,\ldots, k+\ell$) is equal to~$\alpha$ (resp. to
$(y_1+k)(y_2+k)\cdots (y_\ell+k)$).
The set of all shifted shuffles of $\alpha$ and $\beta$ is 
denoted by $\alpha\Cup\beta$.  
The identical permutation $12\cdots k$ is denoted by
$\id_k$.

\proclaim Lemma \ThShuffle. 
On the set $\id_{k_1}\Cup\id_{k_2}\Cup\cdots\Cup\id_{k_r}$ we have 
$$
(\sort\compose\Mc, \El\compose\Mc)\simeq
(\sort\compose\Ic, \El\compose\Ic).
$$ 

{\it Proof}. We construct a bijection 
$\phi: \sigma \mapsto \phi(\sigma)$ on
$\id_{k_1}\Cup\id_{k_2}\Cup\cdots\Cup\id_{k_r}$ satisfying
$\sort\compose\Mc\sigma= \sort\compose\Ic\phi(\sigma)$ and
$\El\compose\Mc\sigma= \El\compose\Ic\phi(\sigma)$.
By induction,
let $\sigma\in\id_{k}\Cup\beta$ and $\Mc(\sigma)=d_1d_2\cdots d_k \Mc(\beta)$
with $\beta\in \id_{k_2}\Cup\cdots\Cup\id_{k_r}$ and $k=k_1$.
As proved in Lemma 6.5 in [\RefHNT], the mapping
$(d_1d_2\cdots d_k,\beta)\mapsto (\sort(d_1d_2\cdots d_k), \beta)$
is bijective. We then define
$\phi(\sigma)=\IcInv\bigl(\sort(d_1d_2\cdots d_k)\Ic(\phi(\beta))\bigr)$. 
We have
$$
\eqalignno{
\sort\compose\Mc(\sigma)&=\sort\bigl(d_1d_2\cdots d_k \Mc(\beta)\bigr) &\cr
&=\sort\bigl(d_1d_2\cdots k_k \sort(\Mc(\beta))\bigr)&\cr
&=\sort\bigl(d_1d_2\cdots k_k \sort(\Ic(\phi(\beta)))\bigr) 
& \hbox{[by induction]}\cr 
&=\sort\bigl(\sort(d_1d_2\cdots k_k) \Ic(\phi(\beta))\bigr) &\cr 
&=\sort\compose\Ic\phi(\sigma). &\cr 
\El\compose\Mc(\sigma)&=\El\bigl(d_1d_2\cdots d_k \Mc(\beta)\bigr) &\cr
&=\El\bigl(\sort(d_1d_2\cdots d_k) \Mc(\beta)\bigr) &\hbox{[by 
Lemma \ThSortMc]}\cr
&=\El\bigl(\sort(d_1d_2\cdots d_k) \Ic(\phi(\beta))\bigr) 
\hphantom{\hbox{[by Lemma \ThElRec]}}
& \hbox{[by Lemma \ThElRec]}\cr
&=\El\compose\Ic\phi(\beta).\qed &\cr
}
$$

{\it Proof of Theorem 1}. 
As used on several occasions (see, e.g., [\RefHNT, Eq. (10)] or [\RefDW, \S 3]),
we have
$$
\eqalign{
&\id_{k_1}\Cup\id_{k_2}\Cup\cdots\Cup\id_{k_r}\cr
&\quad=\{\sigma\mid\ \Iligne(\sigma)\subseteq\{k_1, k_1+k_2, \ldots, 
k_1+k_2+\cdots +k_{r-1} \}\}.\cr
}
$$
By Lemma \ThShuffle\ 
$$
(\sort\compose\Mc, \El\compose\Mc)\simeq
(\sort\compose\Ic, \El\compose\Ic)
$$ 
on the set
$\{\sigma\mid\ \Iligne(\sigma)\subseteq\{k_1, k_1+k_2, \ldots, 
k_1+k_2+\cdots +k_{r-1} \}\}$. It is also true on the set
$\{\sigma\mid\ \Iligne(\sigma)=\{k_1, k_1+k_2, \ldots, 
k_1+k_2+\cdots +k_{r-1} \}\}$ by the inclusion-exclusion principle.
\qed

\section{\SecCompl. The ``Complement" property of the Majcode} 

We rephrase the statement of Theorem 2 as follows.

\proclaim Theorem \ThCompl'.
For each permutation $\sigma$ of $12\cdots n$ let
$$
\tau=\Majcode^{-1}\compose\delta\compose\Majcode(\sigma).\leqno{(R3')}
$$
Then
$$\Ligne\tau=\{1,2,\ldots, n-1\}\setminus\Ligne\sigma.\leqno{(R2)}$$

For example,
take $n=9$ and $\sigma=935721468$. 
Then
$\Ligne\sigma=\{1,4,5\}$,
$\Majcode\sigma=012020203$ and 
$\delta\Majcode\sigma=000325475$.
We have
$\tau=\Majcode^{-1}(000325475)=795128643$.
We verify that
$$\Ligne\tau=\{2,3,6,7,8\}=\{1,2,3,4,5,6,7,8\}\setminus\Ligne\sigma.$$
\medskip

\goodbreak
{\it Proof of Theorem 2'}.
Proceed by induction on the order of the permutation. 
Let $\sigma=x_1x_2\cdots x_n\in\Sym_n$ be a permutation and 
$\sigma'\in\Sym_{n-1}$ be
the permutaion obtained from $\sigma$ by erasing the letter $n$.
Let $\Majcode\sigma=c_1c_2\ldots c_{n-1}c_n$. Then 
$\Majcode\sigma'=c_1c_2\ldots c_{n-1}$. Let $\tau=y_1y_2\cdots y_n\in\Sym_n$ 
be the
permutation defined by relation $(R3')$, 
\  
i.e., 
$\Majcode\tau=d_1d_2\ldots d_{n-1}d_n$ with $d_i=i-1-c_i$ for
$1\leq i \leq n$.
Let $\tau'\in\Sym_{n-1}$ be the permutation obtained from $\tau$ by erasing the
letter $n$. Then
$\Majcode\tau'=d_1d_2\ldots d_{n-1}$.
It is easy to see the $\sigma'$ and $\tau'$  also satisfy relation
$(R3')$. By induction we have
$$
\Ligne\sigma'=\{1,2,\ldots, n-2\} \setminus  \Ligne\tau'.  \leqno{(R2')}
$$

Recall the classical construction of the Majcode consisting of labelling
the slots (see, for example, [\RefRa]). 
Let $\sigma'=x'_1 x'_2 \ldots x'_{n-1}$. Let $x'_0=x'_n=0$
so that the word $x'_0x'_1x'_2\ldots x'_{n-1}x'_n$ has $n$ slots
$(i-1, i)$ with $1\leq i\leq n$.
A slot $(i-1, i)$ is called a {\it descent} (resp. {\it rise}) if
$x'_{i-1} >x'_i$ (resp. $x'_{i-1} <x'_i$).
We label the $k$ descent slots $0,1,2,\ldots,k-1$
from right to left 
and the remaining $n-k$ rise slots 
$k, k+1, \ldots, n-1$
from left to right. 
For $1\leq i\leq n$ let $c_n(i)$ be the label of the slot $(i-1, i)$
and $\sigma^{<i>}\in\Sym_n$ be the permutation obtained from $\sigma'$ 
by inserting $n$ into the slot $(i-1, i)$.
The basic property is that $c_n(i)=\maj\sigma^{<i>} -\maj\sigma'$.
In the same manner, let $d_n(i)$ (for $1\leq i\leq n$) be the label 
of the slots in~$\tau'$.
Thanks to relation $(R2')$
the above construction of labels implies the following simple relation
between $c_n(i)$ and $d_n(i)$~:
$$
c_n(n)=d_n(n)=0 \hbox{\quad and\quad} 
c_n(i)+d_n(i)=n \hbox{\quad (for $ 1\leq i \leq n-1$)}. \leqno{(R4)}
$$
For example, take $\sigma'=35721468$ and $\tau'=75128643$
as in the above example, we have~:
$$
{
\def\tvi{\vrule height 10pt depth 5pt width 0pt}
\def\tv{\tvi\vrule}
\def\m{\kern1pt\nearrow\kern1pt}
\def\d{\kern1pt\searrow\kern1pt}
\vbox{\offinterlineskip\halign{
\tvi\hfil$#$\hfil\ &&\strut\ \hfil$#$\hfil \cr
\hbox{slot of $\sigma'$}&:  &0&\m&3&\m&5&\m&7&\d&2&\d&1&\m&4&\m&6&\m&8&\d&0& \cr
\hbox{label $c_n(i)$}&: & &3 & & 4 & & 5 && 2 && 1 && 6 && 7 && 8 && 0 && \cr
}}
}
$$
$$
{
\def\tvi{\vrule height 10pt depth 5pt width 0pt}
\def\tv{\tvi\vrule}
\def\m{\kern1pt\nearrow\kern1pt}
\def\d{\kern1pt\searrow\kern1pt}
\vbox{\offinterlineskip\halign{
\tvi\hfil$#$\hfil\ &&\strut\ \hfil$#$\hfil \cr
\hbox{slot of $\tau'$}&:&0&\m&7&\d&5&\d&1&\m&2&\m&8&\d&6&\d&4&\d&3&\d&0& \cr
\hbox{label $d_n(i)$}&: & &6 & & 5 & & 4 && 7 && 8 && 3 && 2 && 1 && 0 && \cr
}}
}
$$

If $x_s=n$ and $y_t=n$, that means that the permutation $\sigma$ (resp. $\tau$)
can be constructed by inserting $n$ into the slot
$(s-1,s)$ in $\sigma'$ (resp. slot $(t-1, t)$ in $\tau'$). 
By relation $(R3')$ we have 
$$d_n(t)=n-1-c_n(s).\leqno{(R5)}$$
From $(R4)$ and $(R5)$ we obtain a relation between $c_n(s)$ and
$c_n(t)$~:
$$
c_n(s)=\cases{
n-1, &if $c_n(t)=0$; \cr
c_n(t)-1, &if $c_n(t)\geq  1$.\cr
}\leqno{(R6)}
$$
In fact, relation $(R6)$ gives an algorithm for computing $t$ from $s$.

(st1) if the slot $s$ (that means the slot $(s-1, s)$) is a rise, 
but not the rightmost rise,
then $t$ is the {\it next} rise on the right of $s$. 

(st2) if the slot $s$ is the rightmost rise,  then $t$ is the rightmost slot.

(st3) if the slot $s$ is a descent, 
but not the leftmost descent,
then $t$ is 
first descent preceding $s$ on the left.

(st4) if the slot $s$ is the leftmost descent,  then $t$ is the leftmost slot.

We summarize those cases in the following tableau.

$$
{
\def\tvi{\vrule height 16pt depth 6pt width 0pt}
\def\tv{\tvi\vrule}
\def\m{\kern1pt\nearrow\kern1pt}
\def\d{\kern1pt\searrow\kern1pt}
\vbox{\offinterlineskip\halign{
\tvi\hfil$#$\hfil &&\strut \hfil$\ #$\hfil \cr
\noalign{\hrule}
\tv &&  \tv &&\sigma'&& \tv && \tau' &&\tv \cr
\noalign{\hrule}
\tv &\quad\hbox{(st1)}\quad& \tv 
  & \quad\cdots& {}^s\kern-8pt\m\d\d\d\d{}^t\kern-8pt\m &\cdots \quad 
  & \tv 
  & \quad \cdots& \d\kern-8pt{}^s\m\m\m\m\d\kern-8pt{}^t &\cdots \quad 
  &\tv \cr
\tv &\quad\hbox{(st2)}\quad& \tv 
  & \quad\cdots& {}^s\kern-8pt\m\d\d\d\d\d\kern-8pt{}^t&   \quad
  & \tv 
  & \quad \cdots& \d\kern-8pt{}^s\m\m\m\m\d\kern-8pt{}^t&  \quad
  &\tv \cr
\tv &\quad\hbox{(st3)}\quad& \tv 
  & \quad \cdots& \d\kern-8pt{}^t\m\m\m\m\d\kern-8pt{}^s&\cdots \quad 
  & \tv 
  & \quad\cdots& {}^t\kern-8pt\m\d\d\d\d{}^s\kern-8pt\m&\cdots \quad 
  &\tv \cr
\tv &\quad\hbox{(st4)}\quad& \tv 
  & \quad & {}^t\kern-8pt\m\m\m\m\m\d\kern-8pt{}^s& \cdots \quad
  & \tv 
  & \quad & {}^t\kern-8pt\m\d\d\d\d{}^s\kern-8pt\m& \cdots  \quad
  &\tv \cr
\noalign{\hrule}
}}
}
$$
In each case inserting the letter $n$ into the slot $s$ of $\sigma'$
and inserting the letter $n$ into the slot $t$ of $\tau'$  
produce
two 
permutations $\sigma$ and $\tau$. From the above tableau it easy to see that
 $\sigma$ and $\tau$ satisfy relation $(R2)$. \qed
\medskip

By definitions of ``Invcode", ``Ic", 
``Majcode" and ``Mc" we obtain the following simple
relations between them.
\proclaim Lemma \ThMajcodeMc. 
$$\eqalign{
\Mc&=\bfr\compose\delta\compose\Majcode\compose\bfc,\cr
\Ic&=\bfr\compose\delta\compose\Invcode\compose\bfr\bfi.\cr}
$$

For Example, we have $\Mc (935721468) = 501012010$,  
$\Ic (362715984) = 420520010$, also obtained by the following
calculations.
$$ 
\matrix{
\sigma &=& 935721468 & \cr
\bfc \sigma &=& 175389642 & \cr
\Majcode\compose \bfc \sigma &=& 002135573& \cr
\delta\compose\Majcode\compose \bfc \sigma &=& 010210105& \cr
\bfr \compose \delta\compose\Majcode\compose \bfc \sigma &=& 501012010& \cr
}
$$
$$ 
\matrix{
\sigma &=& 362715984& \cr
\bfi \sigma &=& 531962487& \cr
\bfr \bfi\sigma &=& 784269135& \cr
\Invcode\compose\bfr \bfi\sigma &=& 002320654& \cr
\delta\compose\Invcode\compose\bfr \bfi\sigma &=& 010025024& \cr
\bfr\compose\delta\compose\Invcode\compose\bfr \bfi\sigma &=& 420520010& \cr
}
$$

\medskip
The relation between the statistics ``El" and ``Eul" is given
in the following Lemma.
\proclaim Lemma \ThElEulA. 
Let $d$ be a subdiagonal word of length $n$.
Then 
$$\El(d)=\{1,2,\ldots,n-1\}\setminus\Eul(\delta\bfr(d)).$$

{\it Proof}. We have
$$
\eqalign{
\El(d)&=\Ligne\compose \McInv (d) \cr
&=\Ligne\compose (\bfr\compose\delta\compose\Majcode\compose \bfc)^{-1}(d)\cr
&=\Ligne\compose (\bfc\compose\Majcode^{-1}\compose\delta\compose \bfr)(d)\cr
&=\Ligne (\bfc\compose\Majcode^{-1}(\delta\compose \bfr(d)))\cr
&=\{1,2,\ldots, n-1\}\setminus\Ligne (\Majcode^{-1}(\delta\compose \bfr(d)))\cr
&=\{1,2,\ldots, n-1\}\setminus\Eul(\delta\compose \bfr(d)).\qed\cr
}
$$

In fact, there is another simple, but not trivial, relation between the 
above two statistics.

\proclaim Lemma \ThElEulB. 
Let $d$ be a subdiagonal word of length $n$,
Then 
$$\El(d)=\Eul(\bfr(d)).$$

{\it Proof}.
By Lemma \ThElEulA\ we need verify the following relation 
$$
\Eul(\bfr(d))=\{1,2,\ldots,n-1\}\setminus\Eul(\delta(\bfr(d))).
$$
This is true by Theorem \ThCompl.\qed

\medskip

For every permutation $\sigma$ it is easy to see that
$$
\Invcode\bfr\sigma=\delta\Invcode\sigma.\leqno{(R7)}
$$
\medskip
We end this paper by showing why the equidistribution $(M5)$ 
obtained in [\RefFH] is a special case of Theorem~1. We have
$$
\eqalignno{
\El\compose\Ic\compose\bfi 
&= \Eul\compose\bfr\compose\Ic\compose\bfi &\hbox{[by Lemma \ThElEulB]}\cr
&= \Eul\compose\bfr\bfr\delta\compose\Invcode\compose\bfr\bfi\bfi&
  \hbox{[by Lemma \ThMajcodeMc]}\cr
&= \Eul\compose\delta\compose\Invcode\compose\bfr&\cr
&= \Eul\compose\delta\compose\delta\compose\Invcode&\hbox{[by $(R7)$]}\cr
&= \Eul\compose\Invcode,\cr
}
$$
so that
$$
\eqalignno{
(\Iligne, \Eul\compose\Majcode) 
&\simeq(\Iligne, \Ligne) \cr
&\simeq(\Iligne, \El\compose\Mc) \cr
&\simeq(\Iligne, \El\compose\Ic) \cr
&\simeq(\Ligne, \El\compose\Ic\compose\bfi)\cr
&\simeq(\Ligne, \Eul\compose\Invcode).\cr
}
$$


\vfill\eject
\section{4. Table} 
We give the list of the twenty-four permutations of order 4 and their 
corresponding statistics.
The permutations are sorted according to the statistic ``Iligne".
$$
{
\def\ColEul#1{\ \hbox to 13mm{\hfil#1\hfil}\ \vrule}
\def\ColPerm#1{\ \hbox to 10mm{\hfil#1\hfil}\ }
\def\separ{\noalign{\hrule}}
\vbox{\offinterlineskip\halign{%
\strut \vrule\ColPerm{#}\vrule\ &\ColPerm{#}&\ColEul{#}&\ColPerm{#}&\ColEul{#} 
  &\ColPerm{#}&\ColEul{#}\cr
\separ
$\sigma$&$\Ic$&$\El\compose\Ic$&Mc&$\El\compose\Mc$ & Sc &$\El\compose\Sc$ \cr
\separ
1234&0000&$\epsilon$&0000&$\epsilon$& 0000 &$\epsilon$ \cr
\separ
2134&1000&1&1000&1 & 3000&3 \cr
2314&2000&2&2000&2 & 2000&2 \cr
2341&3000&3&3000&3 & 1000&1 \cr
\separ
1324&0100&1&1100&2 & 0200&2 \cr
1342&0200&2&1200&3 & 0100&1 \cr
3124&1100&2&0100&1 & 2200&13  \cr
3142&1200&3&2200&13 &2100&12  \cr
3412&2200&13&0200&2 &1100&2  \cr
\separ
1243&0010&1&1110&3 & 0010&1 \cr
1423&0110&2&1010&2 & 0110&2 \cr
4123&1110&3&0010&1 & 1110&3 \cr
\separ
3214&2100&12&2100&12 & 3200&23 \cr
3241&3100&13&3100&13 & 1200&3 \cr
3421&3200&23&3200&23 & 3100&13 \cr
\separ
2143&1010&2&2110&13 & 3010&13 \cr
2413&2010&12&0110&2 & 1010&2 \cr
2431&3010&13&3110&23 & 2010&12 \cr
4213&2110&13&2010&12 & 3110&23 \cr
4231&3110&23&3010&13 & 2110&13 \cr
\separ
1432&0210&12&2210&23 & 0210&12 \cr
4132&1210&13&1210&13 & 1210&13 \cr
4312&2210&23&0210&12 & 2210&23 \cr
\separ
4321&3210&123&3210&123 & 3210&123 \cr
\separ
}}
}
$$
From this table we can check the following equidistributions:
$$
\eqalign{
(\Iligne,\sort\compose\Mc, \El\compose\Mc)&\simeq 
(\Iligne,\sort\compose\Ic, \El\compose\Ic),\cr
(\Iligne, \sort\compose\Mc) & \simeq (\Iligne, \sort\compose\Sc),\cr
(\Iligne,\sort\compose\Mc, \El\compose\Mc) &\not\simeq 
(\Iligne,\sort\compose\Sc, \El\compose\Sc).\cr
}
$$


\vfill\eject
\vglue 2mm
\bigskip
{
\eightpoint
\centerline{\bf References} 
\def\thevskip{\smallskip}

\bigskip

\article \RefDW|Jacques D\'esarm\'enien, Michelle L.
Wachs|Descent Classes of Permutations with a Given Number of
Fixed Points|J. Combin. Theory, Ser.~A|64|1993|311--328|
\thevskip

\article \RefFo|D. Foata|On the Netto inversion 
number of a sequence|Proc. Amer. Math. Soc.|19|1968|236--240|
\thevskip

\article \RefFH|D. Foata, G.-N. Han|Une nouvelle transformation pour les 
statistiques Euler-mahoniennes ensemblistes|
Moscow Math. Journal|4|2004|131--152|
\thevskip

\article \RefFS|D. Foata,
M.-P. Sch\"utzenberger|Major Index and Inversion
number of Permutations|Math. Nachr.|83|1978|143--159|
\thevskip

\article \RefGG|Adriano M. Garsia and Ira Gessel|Permutations Statistics and
Partitions|Adv. in Math.|31|1979|288--305|
\thevskip

\article \RefHa|Guo-Niu Han|Distribution Euler-mahonienne: une
correspondance|C. R. Acad. Sci. Paris|310|1990|311--314|
\thevskip

\divers \RefHNT|F. Hivert, J.-C. Novelli, J.-Y. Thibon|Multivariate 
generalizations of the Foata-Sch\"utzenberger equidistribution, 
{\it preprint on arXiv}, {\oldstyle 2006}, 17 pages|
\thevskip

\divers \RefLe|D. H. Lehmer|Teaching combinatorial tricks to a computer, Proc.
Sympos. Appl. Math., vol.~{\bf 10}, {\oldstyle1960}, p. 179--193. Amer. Math. 
Soc., Providence, R.I|
\thevskip

\livre \RefLoI|M. Lothaire|Combinatorics on Words|Addison-Wesley,
London {\oldstyle 1983} (Encyclopedia of Math. and its Appl., {\bf
17})|
\thevskip

\livre \RefLoII|M. Lothaire|Algebraic Combinatorics on Words|Cambridge Univ.
Press, {\oldstyle 2002} (Encyclopedia of Math. and its Appl., {\bf 90})|
\thevskip

\livre \RefMac|P. A. MacMahon|Combinatory Analysis, {\rm vol.~1}|Cambridge Univ.
Press, {\oldstyle 1915}|
\thevskip

\article \RefRa|Don P. Rawlings|Generalized Worpitzky identities with
applications to permutation enumeration|Europ. J. 
Comb.|2|1981|67-78|
\thevskip

\divers \RefSk|Mark Skandera|An Eulerian partner for inversions, {\sl S\'em.
Lothar. Combin.}, vol.~{\bf 46}, [B46d], {\oldstyle 2001}, 19 pages.
{\tt http://www.mat.univie.ac.at/\char126slc}|
\thevskip

} 

\vfill\eject

\bye